\documentclass[12pt]{article}
\usepackage[cp1251]{inputenc}
\usepackage{amssymb}
\usepackage{amsthm}
\usepackage{inputenc}
\begin{document}

\begin{center}
\textbf{ORLICZ SPACES ASSOCIATED WITH A SEMI-FINITE VON NEUMANN ALGEBRA}
\end{center}

\vspace{0.5 cm}
\begin{center}
\textbf{Sh. A. Ayupov$^{1},$ V.I. Chilin$^2,$ R. Z. Abdullaev$^3$}
\end{center}

$^1$ \emph{Institute of Mathematics and Information Technologies, Uzbekistan Academy of Science,
Dormon yoli, 29, 100125, Tashkent. Uzbekistan}

and

\emph{The Abdus Salam International Centre for Theoretical Physics. Triest. Italy.}

Coresponding author. e-mail: \emph{sh\_ayupov@mail.ru}

 $^{2}$ \emph{National University of  Uzbekistan,} e-mail: \emph{chilin@usd.uz}

$^3$\emph{ Tashkent State Pedagogical University,} e-mail: \emph{arustambay@yandex.ru}

\begin{center}
\textbf{Introduction}
\end{center}

Construction and investigation of various
classes of symmetric spaces of measurable operators affiliated with a von Neumann
algebra $M$ is one of important applications of the non commutative integration theory for a faithful
normal semi-finite trace  on the von Neumann algebra $M.$ Examples of such spaces are given,
in particular by non commutative $L_p$-spaces $L_p(m,\tau)$ \cite{Yea2} and by Orlicz spaces
$L_\Phi(M,\tau)$ associated with an $N$-function $\Phi$
\cite{Kun},\cite{Mur1},\cite{Mur2}. All these spaces are realized as ideal subspaces of the $*-$algebra
$S(M)$ of measurable operators affiliated with $M.$

Investigations based on the modular theory for von Neumann algebras enable to consider non commutative versions of
$L_p$-spaces associated with states and weights (see e.g. the survey \cite{Tru4}).But in these cases in general
$L_p$-spaces can not be realizes as ideal subspaces of $S(M).$ This fact explains in particular why in their attempt
to introduce non commutative Orlicz spaces for states in \cite{Ras} as a subspaces of $S(M),$ the authors we unable
to prove the completeness of these spaces with respect to the Luxemburg norm.

In the present paper we introduce a certain class of non commutative Orlicz spaces,
associated with arbitrary faithful normal locally-finite weights on a semi-finite von Neumann algebra $M.$
We describe the dual spaces for such Orlicz spaces and, in the case of regular weights , we show that
they can be realized as linear subspaces of the algebra of $LS(M)$ of locally measurable operators affiliated with $M.$

For the terminology and notations from the von Neumann algebras theory we refer to \cite{Tak} and from theory of measurable and locally
measurable operators refer to \cite{Tak},\cite{Yea1}.

\begin{center}
\textbf{Preliminaries}
\end{center}

Let $M$ be  a von Neumann algebra
acting on a Hilbert space $H$ with $\mathbf{1}$-the identity operator on $H,$ and let
$P(M)=\{p\in M: p=p^2=p^*\}$ be the lattice of all projection from $M$. Denote by $S(M)$ (respectively by $LS(M)$)
the $^*$-algebra of all measurable (respectively, locally-measurable) operators affiliated with $M$.
It is well-known that $S(M)$ is a $^*$-subalgebra in $LS(M)$, and $M$ is a $^*$-subalgebra of $S(M)$
(\cite{Chl}, Ch.2).

If $x\in LS(M)$ and $x=u|x|$ is its polar decomposition, where
$|x|=(x^*x)^{1/2}$ and $u$ is a partial isometry, then we have that $u\in M$ and $|x|\in LS(M).$
It is also known that the spectral family of projections
$e_{\lambda}(x)$ for a self adjoint operator  $x\in LS(M)$, always belongs to
$P(M).$

Given a subset
$A\subset LS(M)$, put

 $A_h=\{x\in A: x=x^*\},$ and

 $A_+=\{x\in A: (x\xi,\xi)\geq0\}$ for all $\xi\in D(x),$

 where $D(x)$ is
the domain of the operator  $x\in LS(M),$ and $(\cdot,\cdot)$ is the inner product in the Hilbert space $H.$

Let  $\tau$ be a faithful normal semi-finite trace on  $M.$
For each real number $p\geq1$ consider the set
$$
L_p(M,\tau)=\left\{x\in S(M): \int\limits_0^\infty\lambda^pd\tau(e_\lambda(|x|))<\infty\right\}.
$$
It known \cite{Yea2} that  $L_p(M,\tau)$ is a linear subspace in $S(M)$ and the function
$\|x\|_p=\left(\int\limits_0^\infty\lambda^pd(\tau (e_\lambda(|x|))\right)^{1/p}$
is a norm, which turns  $L_p(M,\tau)$ into a Banach space.

A map $\varphi : M_+\rightarrow[0,\infty]$ is said to be \emph{a weight} if
$$
\varphi(x+y)=\varphi(x)+\varphi(y), \ \ \varphi(\lambda
x)=\lambda\varphi(x), \ \ (x,y\in M_+, \lambda\geq0, \textrm{where, }
0\cdot\infty=0).
$$

A weight $\varphi$ is said to be

--- \emph{normal}, if $\varphi(x)=\sup\varphi(x_i) (x_i\nearrow x; x_i,x\in M_+);$

--- \emph{faithful}, if $\varphi(x)=0,$ $x\in M_+$ implies that $x=0;$

---\emph{ semi-finite}, if the linear span  $m_\varphi$ of the cone  $m_\varphi^+=\{x\in M_+ : \varphi(x)<\infty\}$
is dense in  $M$ with respect to the ultra-weak topology;

--- \emph{locally finite}, if
$$
\forall x\in M_+ \ \ (x\neq0) \ \ \exists y\in M_+:y\leq x, 0<\varphi(y)<\infty;
$$

---\emph{ regular}, if
$$
\forall\omega\in (M_*)_+\,\, (\omega\neq0) \ \exists\omega'\in (M_*)_+\,\,
(\omega'\neq0): \omega'\leq\omega, \omega'\leq\varphi,
$$
where $(M_*)_+$ is the set of all positive ultra-weakly continuous linear functionals on $M.$

If the weight $\varphi$ is a trace, i.e. when  $\varphi(x^*x)=\varphi(xx^*)$ for all $x\in M,$
the properties of semi-finiteness and locally finiteness (and respectively of faithfulness and regularity)
of $\varphi$ coincide with each other \cite{Tru3}.

For a faithful normal semi-finite weight  $\varphi$ on $M$ there exists a
uniquely defined non singular self-adjoint positive operator $h$, affiliated
with $M$ such that $\varphi(\cdot)=\tau(h\cdot )$ , and which is called the Radon-Nikodym derivative
of the weight  $\varphi$ with respect is the trace $\tau$ \cite{Ped}.

Recall the following result

\textbf{Theorem 1.}\cite{Tru3}
\emph{Let  $\tau$ be a faithful normal semi-finite trace on  $M$ and let
$\varphi=\tau(h\cdot)$  be a faithful normal semi-finite weight on $M,$ where $h$ is the Radon-Nikodym derivative
of $\varphi$ with respect to $\tau.$
Then}

$(i)$ \emph{the weight $\varphi$ is locally finite if and only if the operator $h$ is locally measurable;}

$(ii)$ \emph{the weight $\varphi$ is regular if and only if the operator
 $h^{-1}$ is locally measurable.}

Now let $\varphi(\cdot)=\tau(h\cdot)$ be a faithful  normal locally finite weight on $M.$
For real numbers $p\geq1$ and $\alpha\in[0,1]$ put
$$
m_\alpha^{1/p}=\{x\in M: h^{\alpha/p}xh^{(1-\alpha)/p}\in L_p(M,\tau)\};
$$
$$
\|x\|_{p,\alpha}=\|h^{\alpha/p}xh^{(1-p)/p}\|_{p}.
$$

In  \cite{Tru2} it has been proved that  $m_\alpha^{1/p}$ is a linear subspace in $M,$ and $\|\cdot\|_{p,\alpha}$
is a norm on  $m_\alpha^{1/p}.$
The completion of the normed space $(m_\alpha^{1/p},\|\cdot\|_{p,\alpha})$ is denoted by  $L_p(M,\varphi).$
In  \cite{Tru2} it is proved that the Banach space $(L_p(M,\varphi),\|\cdot\|_{p,\alpha})$
is isometrically isomorphic to the space
$(L_p(M,\tau),\|\cdot\|_{p})$ for all $\alpha\in[0,1].$

In order to define the Orlicz space associated with a weight, we need the notion of $N$-function.

A continuous non-negative convex monotone increasing function  $\Phi$ on the set of real numbers $\mathbb{R}$
is called $N$-\emph{function} \cite{Kra}, if
$$
\Phi(t)=\int\limits_{0}^{|t|}p(s)ds,
$$
where $p(s)$ is a non-decreasing function, positive for $s>0$ and right continuous for
$s\geq0$ , which satisfies the conditions
$$
p(0)=0, \ \
p(\infty)=\lim\limits_{s\rightarrow\infty}p(s)=\infty.
$$

For each $N$-function $\Phi(t)$  a complementary $N$-function
 $\Psi(t)$
is defined as
$$
\Psi(t)=\int\limits_{0}^{|t|}q(s)ds,
$$
where $q(s)=\sup\{t\geq0: p(t)\leq s\}.$
It is clear that the complementary $N$-function for the $N$-function
$\Psi(t)$ coincides with the initial function $\Phi(t),$
and moreover the following Young inequality is valid
$$
ts\leq\Phi(t)+\Psi(s) \ \ \textrm{for all}\ \  t,s\geq0.
$$

We say that an $N$-function $\Phi(t)$
satisfies the  $(\delta_2,\Delta_2)$-condition, if given any real $k>0$
there exists a positive number $r(k)$  such that $\Phi(kt)\leq
r(k)\Phi(t)$ for all $t\geq 0.$
Examples of $N$-function which satisfy the $(\delta_2,\Delta_2)$-condition
are given by the function $\Phi(t)=\frac{1}{p}|t|^p,$ $p>1.$

Let $\Phi(t)$ be an $N$-function and let $x\in LS_h(M),$ $x=\int\limits_{-\infty}^\infty\lambda de_\lambda(x).$
It is known (\cite{Chl}, \S 2.3) that one can define a self-adjoint operator
$\Phi(x)=\int\limits_{-\infty}^\infty\Phi(\lambda) de_\lambda(x),$
and moreover $\Phi(x)\in LS(M).$

Let us extend the faithful normal semi-finite trace $\tau$ from $M_+$ to operators from  $LS_+(M)$
as
$$
\tau(x)=\sup\limits_{t\geq1}\tau\left(\int\limits_{0}^t\lambda de_\lambda(x)\right)=
\int\limits_{0}^\infty\lambda d\tau (e_\lambda(x)).
$$
It is known (e.g. \cite{Chl}, \S 4.1), that
$$
\tau(x)=\sup\{\tau(y): y\in M_+, y\leq x\}
$$
for all $x\in LS_+(M).$

It is clear that $\tau(|x|)<\infty$ for $x\in LS(M)$ if and only if $x\in L_1(M,\tau);$
in this case $\tau(\mathbf{1}-e_\lambda(|x|))<\infty$ for all $\lambda>0.$
Further we shall need the following result.

\textbf{Proposition 1.} \cite{Fac} \emph{If $x,y\in LS_+(M),$ then}

$(i)$ $\tau(f(x))\leq\tau(f(y))$ for $x\leq y$ \emph{for each continuous monotone increasing
function $f:[0,\infty)\rightarrow \mathbb{R}$
with $f(0)=0;$}

$(ii)$ $\tau(f(\lambda x+(1-\lambda)y))\leq\lambda\tau(f(x))+(1-\lambda)\tau(f(y))$ \emph{for all
$\lambda\in [0,1]$ and each convex monotone increasing function $f$ with $f(0)=0.$}

Let $\Phi$ be an $N$-function. The set

$K_\Phi=\{x\in S(M): \tau(\Phi(|x|))\leq1\}$

is an absolutely convex subset in $S(M)$ \cite{Kun}.
The linear subspace $L_\Phi(M,\tau)=\bigcup\limits_{n=1}^\infty nK_\Phi$ equipped with the norm
$$
\|x\|_\Phi=\inf\left\{\lambda>0: \frac{x}{\lambda}\in K_\Phi\right\}, \eqno(1)
$$
is a Banach space \cite{Kun} which is called the Orlicz space associated with $M, \tau$ and $\Phi$.
If the $N$-function $\Phi$ satisfies the $(\delta_2,\Delta_2)$-condition, then
$$
L_\Phi(M,\tau)=\{x\in S(M): \tau(\Phi(|x|))<\infty\},
$$
moreover the linear subspace

$m_\Phi^\tau=\{x\in M: \tau(\Phi(|x|))<\infty\}$

is dense in
 $\left(L_\Phi(M,\tau), \|\cdot\|_\Phi\right).$

Note that
$$
m_\tau=\{x\in M: \tau(|x|)<\infty\}\subset m_\Phi^\tau. \eqno(2)
$$

Indeed, from the equalities
$$
\lim\limits_{t\downarrow 0}\frac{\Phi(t)}{t}=\lim\limits_{t\downarrow 0}p(t)=0
$$
it follows that $\Phi(t)\leq t$ for sufficiently small $t>0.$

Therefore for $x\in m_\tau$ there exists $t_0>0$ such that

$\tau(\Phi(|x|e_{t_0}(|x|)))=\int\limits_0^{t_0}\Phi(\lambda)d\tau(e_\lambda(|x|))\leq
\int\limits_0^{t_0}\lambda d\tau(e_\lambda(|x|))=\tau(|x|e_{t_0}(|x|))<\infty.$

Since $\tau(\mathbf{1}-e_{t_0}(|x|))<\infty,$ we have that

$\tau(\Phi(|x|(\mathbf{1}-e_{t_0})))\leq
\Phi(\|x\|_M)\tau(\mathbf{1}-e_{t_0}(|x|))<\infty,$

where $\|\cdot\|_M$ is the $C^*$-norm on $M.$
Therefore $\tau(\Phi(|x|))<\infty,$ i.e. $x\in m_\Phi^\tau.$

\textbf{Proposition 2.}
\emph{If the $N$-function $\Phi$ satisfies the $(\delta_2,\Delta_2)$-condition, then $m_\tau$
is dense in }$L_\Phi(M,\tau).$

\emph{Proof.}
Since
 $m_\tau\subset m_\Phi^\tau$ (see (2)) and $m_\Phi^\tau$ is dense in
$L_\Phi(M,\tau),$ it sufficient to prove that $m_\tau$ is dense in $m_\Phi^\tau.$
Moreover since each element of
$m_\phi^\tau$ is a finite linear combination of positive elements from $m_\Phi^\tau$
it sufficient to show that every element  from $x\in\left(m_\Phi^\tau\right)_+$
belongs to the closure  of $m_\tau$ in $L_\Phi(M,\tau).$ First, let us show that
$$
x_n=x(\textbf{1}-e_\frac{1}{n})\in m_\tau,
$$
where  $e_\lambda=e_\lambda(x),$ $\lambda>0$ is the spectral family of $x.$
From
$$
\Phi\left(\frac{1}{n}\right)\tau\left(\textbf{1}-e_{\frac{1}{n}}\right)=\tau\left(\Phi\left(\frac{1}{n}
\left(\textbf{1}-e_{\frac{1}{n}}\right)\right)\right)\leq
\tau\left(\Phi\left(x\left(\textbf{1}-e_{\frac{1}{n}}\right)\right)\right)\leq \tau(\Phi(x))<\infty,
 $$
it follows that  $\tau\left(\textbf{1}-e_\frac{1}{n}\right)<\infty$ and the inequality
$0\leq x\left(\textbf{1}-e_\frac{1}{n}\right)\leq
\|x\|_M\left(\textbf{1}-e_\frac{1}{n}\right)$ implies that
$$
x_n=x\left(\textbf{1}-e_\frac{1}{n}\right)\in m_\tau.
$$

Since $0\leq xe_{\frac{1}{n}}\downarrow0$ when $n\rightarrow\infty,$ it follows that
$\tau\left(\Phi\left(\frac{1}{\varepsilon}xe_{\frac{1}{n}}\right)\right)\downarrow0$ for any $ \varepsilon>0.$
In particular, there exists $n(\varepsilon)$ such that
$\tau\left(\Phi\left(\frac{1}{\varepsilon}xe_{\frac{1}{n}}\right)\right)<1$
for $n\geq n(\varepsilon),$ i.e. $\left\|xe_{\frac{1}{n}}\right\|_\Phi<\varepsilon.$ This means that
 $\|x-x_n\|_\Phi\rightarrow0$, i.e. $m_\tau$ is dense $m_\Phi^\tau.$

 The proof is complete.\ $\Box$\\[-2mm]

Let $\Psi$ be the complementary $N$-function for the  $N$-function $\Phi$ satisfying the  $(\delta_2,\Delta_2)$-condition.
In this case given any $y\in L_\Psi(M,\tau)$ the function $f_y(x)=\tau(xy), x\in L_\Phi(M,\tau),$
defines the general form of continuous linear functionals on  $L_\Phi(M,\tau)$ \cite{Kun}, moreover
$$
\|f_y\|=\sup\{|\tau(xy)|: x\in L_\Phi(M,\tau), \|x\|_\Phi\leq1\}=\|y\|_\Psi.
$$
Further we shall need also two inequalities from the following proposition.

\textbf{Proposition 3.}
\emph{Let $\tau$ be a faithful normal semi-finite trace on a von Neumann algebra $M$. Then}

$(i)$(\cite{Chl}, \S 3.4).
\emph{ Given any $x,y\in LS(M)$ there exist  two partial isometries $u,v\in M$ such that}
$$|x+y|\leq u^*|x|u+v^*|y|v.$$

$(ii)$ \cite{Bra}. \emph{For every  $N$-function  $\Phi$, arbitrary operator $z\in M$ with $\|z\|_M\leq1$, and for each
 $x\in LS_+(M)$ we have the following inequality}
$$
\tau(\Phi(z^*xz))\leq\tau(z^*\Phi(x)z).
$$

\begin{center}
\textbf{Orlicz spaces associated with a weight}
\end{center}

In this section an approach is suggested for the construction of Orlicz spaces associated
with a faithful normal locally finite weight on a semi-finite von Neumann algebra for an
$N$-function satisfying the $(\delta_2,\Delta_2)$-condition. For these spaces the dual spaces are described. In the
case of regular locally finite normal weights the constructed Orlicz spaces are represented as spaces
of locally measurable operators.

Let $\tau$ be a faithful normal semi-finite trace on a von Neumann algebra $M$. From now on $\varphi$
denotes a faithful normal locally finite weight on $M$. Therefore the Radon-Nikodym derivative
$h$ of the weight $\varphi$ with respect to $\tau$ is a positive locally measurable non-singular operator.

Given an $N$-function $\Phi$ and a real number $\alpha\in[0,1]$ put
$$
U(x)=U_{\Phi,\alpha}^{\varphi,\tau}(x) = (\Phi^{-1}(h))^\alpha x
(\Phi^{-1}(h))^{1-\alpha}, \ \ x\in LS(M).
$$
It is clear that
$U(x)\in LS(M)$ and $\Phi(|U(x)|)\in LS(M).$

Consider the functional on $LS(M)$ defined by
$$
O_{\Phi,\alpha}^{\varphi,\tau}(x)=\tau(\Phi(|U(x)|)),
$$
and put
$$
m_{\Phi,\alpha}^{\varphi,\tau}=\left\{x\in M: O_{\Phi,\alpha}^{\varphi,\tau}(x)<\infty\right\}.
$$
Consider on the set $m_{\Phi,\alpha}^{\varphi,\tau}$ the functional
$$
\|x\|_{\Phi,\alpha}^{\varphi,\tau}=\inf\left\{\lambda>0:
O_{\Phi,\alpha}^{\varphi,\tau}\left(\frac{x}{\lambda}\right)\leq1\right\}.
$$

\textbf{Theorem 2.}
\emph{If the $N$-function $\Phi$ satisfies the
$(\delta_2,\Delta_2)$-condition, then the set $m_{\Phi,\alpha}^{\varphi,\tau}$
is a linear subspace in $M$.}

In order to prove this theorem we need the following inequality.

\textbf{Lemma 1.}
\emph{For the $N$-function
$\Phi$ and real number $\lambda\in[0,1]$
the following inequality is valid}

$$\tau(\Phi(|U(\lambda x)|))\leq \lambda \tau(\Phi(|U(x)|))\eqno(3) $$
\emph{for all} $x\in m_{\Phi,\alpha}^{\Phi,\tau}.$

\emph{Proof.}
By the linearity of the map $U$ we have

$\tau(\Phi(|U(\lambda x)|))=
\tau(\Phi(\lambda|U(x)|)).$

From the inequality  $(ii)$ in Proposition 1 with $y=0$,
we obtain
$$
\tau(\Phi(|U(\lambda x)|))\leq \lambda
\tau(\Phi(|U(x)|)).
$$

The proof of lemma is complete.\ $\Box$\\[-2mm]

\emph{Proof of the theorem 2.}
The inequality (3) above implies that  $m_{\Phi,\alpha}^{\varphi,\tau}$
is closed under the multiplication by complex number  $k$ with $|k|\leq1.$ Let us show that for
  $x\in m_{\Phi,\alpha}^{\varphi,\tau}$ and any complex number $k$ with $|k|>1$ we have
  that $kx\in m_{\Phi,\alpha}^{\varphi,\tau},$ i.e.
$\tau(\Phi(|U(kx)|))<\infty.$

Since $\Phi$ satisfies the
$(\delta_2,\Delta_2)$-condition, given any positive number  $|k|$
there exists a positive number $r(|k|)$  such that
$\Phi(|k|t)\leq r(|k|)\Phi(t)$ for all $t\geq 0$.  Therefore
$(\delta_2,\Delta_2)$-condition implies that
$$
\tau(\Phi(|U(kx)|))=\tau(\Phi(|k||U(x)|))\leq
$$
$$\leq r(|k|) \tau(\Phi(|U(x)|))<\infty,
$$
 i.e. the set
$m_{\Phi,\alpha}^{\varphi,\tau}$ is closed under multiplication by any complex number.

Now let us prove that the sum of any two operators from
$m_{\Phi,\alpha}^{\varphi,\tau}$ also belongs to
$m_{\Phi,\alpha}^{\varphi,\tau}$.
Let $x,y\in m_{\Phi,\alpha}^{\varphi,\tau},$ i.e.
$\tau(\Phi(|U(x)|))<\infty$ and $\tau(\Phi(|U(y)|))<\infty.$
The inequalities $(i)$ and $(ii)$ from proposition 3, the linearity of the operator
$U$, the convexity of  $\Phi$,  the tracial property of $\tau$ and the fact that $m_{\Phi,\alpha}^{\varphi,\tau}$
is closed under the multiplication by complex numbers imply:
$$
\tau(\Phi(|U(x+y)|))=\tau(\Phi(|U(x)+U(y)|))\leq
\tau(\Phi(u^*|U(x)|u+v^*|U(y)|v))\leq
$$
$$
\leq
\tau\left(\frac{1}{2}(\Phi(2u^*|U(x)|u)+\Phi(2v^*|U(y)|v))\right)=
$$
$$
=\frac{1}{2}(\tau\left(\Phi(u^*2|U(x)|u))+\tau(\Phi(v^*2|U(y)|v)\right)\leq
\frac{1}{2}\left(\tau(u^*\Phi(2|U(x)|)u)+\tau(v^*\Phi(2|U(y)|)v)\right)\leq
$$
$$
\leq\frac{1}{2}(\tau(\Phi(2|U(x)|))+\tau(\Phi(2|U(y)|)))=\frac{1}{2}(\tau(\Phi|U(2x)|)+\tau(\Phi|U(2y)|))<\infty,
$$
i.e. $x+y\in m_{\Phi,\alpha}^{\varphi,\tau}.$

The proof is complete.\ $\Box$\\[-2mm]

\textbf{Theorem 3.}
\emph{The set
$$
K_{\Phi,\alpha}^{\varphi,\tau}=\{x\in M: O_{\Phi,\alpha}^{\varphi,
\tau}(x)\leq1\}
$$
is absolutely convex and absorbing in  $m_{\Phi,\alpha}^{\varphi,\tau}.$}

\emph{Proof.} Let us prove the convexity of
$K_{\Phi,\alpha}^{\varphi,\tau}$. Let  $x, y\in
K_{\Phi,\alpha}^{\varphi,\tau} $ and $\lambda\in[0,1].$ In view of Proposition 3$(i)$
there exist partial isometries  $u$ and $v$ in $M$ such that
$$
|\lambda U(x)+(1-\lambda)U(y)|\leq \lambda u^*|U(x)|u+(1-\lambda)v^*|U(y)|v.
$$
From the inequalities of Proposition 1 and 3 and from the tracial property  of  $\tau$ we obtain
$$
\tau(\Phi |\lambda U(x)+(1-\lambda)U(y)|)\leq \lambda \tau( \Phi
(u^*|U(x)|u))+(1-\lambda)\tau(\Phi(v^*|U(y)|v)) \leq
$$
$$
\leq \lambda\tau( u^*\Phi (|U(x)|u))+(1-\lambda)\tau(v^*\Phi(|U(y)|)v)  \leq
 \lambda \tau( \Phi (|U(x)|))+(1-\lambda)\tau(\Phi(|U(y)|)),
$$
i.e.
$$
O_{\Phi,\alpha}^{\varphi,\tau}(\lambda x+(1-\lambda)y)\leq
\lambda O_{\Phi,\alpha}^{\varphi,\tau}(x)+(1-\lambda)O_{\Phi,\alpha}^{\varphi,
\tau}(y),
$$
which implies the convexity of  $K_{\Phi,\alpha}^{\varphi,\tau}$.

The inequality (3) shows that the set  $K_{\Phi,\alpha}^{\varphi,\tau}$
is balanced, and hence is absolutely convex.

Finally let us move that $K_{\Phi,\alpha}^{\varphi,\tau}$ is absorbing in $m_{\Phi,\alpha}^{\varphi,\tau}$.

If  $x\in m_{\Phi, \alpha}^{\varphi, \tau},$ then there exists $ t>1$ such that $O_{\Phi, \alpha}^{\varphi,
\tau}(x)<t.$ Let $\lambda\in\mathbb{C}$ and $|\lambda|\geq t$ By Lemma 1 we have that

$O_{\Phi,
\alpha}^{\varphi, \tau}(\frac{x}{\lambda})\leq\frac{1}{|\lambda|} O_{\Phi,
\alpha}^{\varphi, \tau}(x)\leq\frac{1}{t}O_{\Phi,
\alpha}^{\varphi, \tau}(x)< 1,$

i.e. $\frac{x}{\lambda}\in K_{\Phi,
\alpha}^{\varphi, \tau}(x).$
The proof is complete.\ $\Box$\\[-2mm]

\textbf{Corollary 1.} \emph{The Minkovsky functional of the set
$K_{\Phi,\alpha}^{\varphi,\tau}$  defined as
$$
\|x\|_{\Phi,\alpha}^{\varphi,\tau}=\inf\left\{\lambda>0:
\frac{x}{\lambda}\in K_{\Phi,\alpha}^{\varphi,\tau}\right\},\eqno(4)
$$
is a norm on the linear space $m_{\Phi,\alpha}^{\varphi,\tau}$.}

\emph{Proof.}
It is sufficient to prove that
$\|x\|_{\Phi,\alpha}^{\varphi,\tau}=0$
implies that $x=0$. Indeed, if
$\|x\|_{\Phi,\alpha}^{\varphi,\tau}=0$ then
$O_{\Phi,\alpha}^{\varphi,\tau}\left(\frac{x}{\lambda}\right)\leq1$
 for all $\lambda\in(0,1)$. By Lemma 1 we obtain that
$\frac{1}{\lambda}O_{\Phi,\alpha}^{\varphi,\tau}\left(x\right)
\leq
O_{\Phi,\alpha}^{\varphi,\tau}\left(\frac{x}{\lambda}\right)\leq1$
for all $\lambda\in(0,1),$ i.e.
$O_{\Phi,\alpha}^{\varphi,\tau}(x)=0.$ Faithfulness of
$\tau$ then implies that $\Phi^{-1}(h)^\alpha x \Phi^{-1}(h)^{1-\alpha}=0$.
 Since $h\in LS_+(M)$ (see theorem 1$(i)$) and $h$ is a non singular operator,
we have that $\Phi^{-1}(h)^\alpha,$ $\Phi^{-1}(h)^{1-\alpha}\in LS_+(M)$ and
$\Phi^{-1}(h)^\alpha,$ $\Phi^{-1}(h)^{1-\alpha}$ are non singular operators too.
Let

$\Phi^{-1}(h)=\int\limits_0^\infty\lambda de_\lambda,$
$x_n=\int\limits_{\frac{1}{n}}^n\left(\frac{1}{\lambda}\right)^\alpha\lambda de_\lambda,$
$y_n=\int\limits_{\frac{1}{n}}^n\left(\frac{1}{\lambda}\right)^{1-\alpha}\lambda de_\lambda.$

Using $x_n\Phi^{-1}(h)^\alpha=e_n-e_{\frac{1}{n}}=\Phi^{-1}(h)^{1-\alpha}y_n$ we see that
$\left(e_n-e_{\frac{1}{n}}\right)x\left(e_n-\frac{1}{n}\right)=0.$
Since $\Phi^{-1}(h)$ is a non singular operator, it follows that $\left(e_n-e_{\frac{1}{n}}\right)\uparrow \mathbf{1}$
 for all $n=1,2,\ldots$
when $n\rightarrow \infty.$ Consequently, $x=0.$
The proof is complete.\ $\Box$\\[-2mm]

Denote by
$L_{\Phi,\alpha}(M,\varphi,\tau)$ the Banach space obtained as the completion of
 $m_{\Phi,\alpha}^{\varphi,\tau}$ in the norm
$\|\cdot\|_{\Phi, \alpha}^{\varphi,\tau}$ and call this completion \emph{the Orlicz space}
constructed by the
$N$-function $\Phi$ on the von Neumann algebra $M$ with respect to the faithful normal locally finite weight $\varphi$.
It is clear that if $\varphi$ is a trace or $M$ is a commutative von Neumann algebra, then the norm
 $\|\cdot\|_{\Phi, \alpha}^{\varphi,\tau}$ and the space  $L_{\Phi,\alpha}(M,\varphi,\tau)$ do not depend
 on $\alpha\in[0,1].$

Note also that in the case where  $\Phi(t)=\frac{1}{t}|t|^p,  \ p> 1,$
the norm $\|\cdot\|_{\Phi, \alpha}^{\varphi,\tau}$ and the space
$L_{\Phi,\alpha}(M,\varphi,\tau)$ do not depend
 on the choice of the faithful normal semi-finite trace  $\tau$ and of $\alpha\in[0,1]$ \cite{Tru2}.

For general $N$-functions $\Phi$ this is not true even in the commutative case.

\textbf{Example.}
Take  $M=l_\infty,$ $f_i=\{0,...,0,1,0...\},$ where 1 is on the $i$-th position, and put
$\Phi(t)=|t|^\beta(\ln |t|+1), \ t\neq0, \ \beta>1, \ \Phi(0)=0.$
In (\cite{Kra}, Ch. I, \S 4) it is proved that $\Phi$ is an $N$-function  satisfying the $(\delta_2,\Delta_2)$-condition.
Consider the trace $\nu$ on $l_\infty$
defined as
$\nu(f_i)=\frac{1}{i^2}((e^{i^2})^{2\beta}(2i^2+1))^{-1}.$ Put

$h=\{e^{\beta i^2}+(i^2+1)\}_{i=1}^\infty,$
$f=\Phi^{-1}(h)=\{e^{i^2}\}_{i=1}^\infty.$

Now define the trace $\mu$ on $l_\infty$ as $\mu(\cdot)=\nu(h\cdot).$

Let us show that in this case the norms $\|\cdot\|_{\Phi, 1}^{\mu,\nu}$ and
$\|\cdot\|_{\Phi, 1}^{\mu,\mu}$ are not equivalent on the ideal $E$ of all finite
sequences from $l_\infty$ (it is clear that $E\subset m_{\Phi,\alpha}^{\mu,\nu}$
and $E\subset m_{\Phi,\alpha}^{\mu,\mu}$).
For this it is sufficient to find a sequence $\{x_n\}$
 of elements from $(K_{\Phi,1}^{\mu,\nu})\cap E$ such that
$\{x_n\}\subset\!\!\!\!\!/\lambda K_{\Phi,1}^{\mu,\mu}$
for all
$\lambda>0.$
Let $x_n=\sum\limits_{i=2}^n e^{i^2}f_i.$ It is clear that for
commutative algebras one has
$$
O_{\Phi,1}^{\mu,\nu}(x)=\nu(\Phi(|\Phi^{-1}(h)x|))
$$
and
$$
O_{\Phi,1}^{\mu,\mu}(x)=\mu(\Phi(|x|)).\eqno(5)
$$
Therefore
$$
O_{\Phi,1}^{\mu,\nu}(x_nf_i)=\nu(\Phi(fx_nf_i))=\nu(\Phi((e^{2i^2}f_i)^2))=(e^{i^2})^{2\beta}(2i^2+1)\nu(f_i)=\frac{1}{i^2}.
$$
Hence
$$
O_{\Phi,1}^{\mu,\nu}(x_n)=\sum\limits_{i=2}^n \frac{1}{i^2}<1,
\textrm {i.e.}\,  x_n\in K_{\Phi,1}^{\mu,\nu}, \eqno(6)
$$
for all $n.$
Let us show that $\{x_n\}\subset\!\!\!\!\!/\lambda K_{\Phi, 1}^{\mu, \mu}$
for all positive real $\lambda.$ From (5) we have
$$
O_{\Phi,1}^{\mu,\mu}(x_nf_i)=\mu(\Phi(x_nf_i))=
\nu(h\Phi(x_n)f_i)=\nu(\Phi(\Phi^{-1}(h))\Phi(x_n)f_i)=\nu(\Phi(f)\Phi(x_n)f_i)=
$$
$$
=(e^{i^2})^{2\beta}(i^4+2i^2+1)\nu(f_i)>(e^{i^2})^{2\beta}i(2i^2+1)\nu(f_i)=\frac{1}{i}.
$$
Therefore  $O_{\Phi,1}^{\mu,\mu}(x_n)>\sum\limits_{i=2}^n\frac{1}{n},$
and hence
$$
\{x_n\}\subset\!\!\!\!\!/\lambda
K_{\Phi,1}^{\mu,\mu}\eqno(7)
$$
for all positive
$\lambda.$ From (6) and (7) it follows that the norms
$\|\cdot\|_{\Phi,1}^{\mu,\nu}$ and $\|\cdot\|_{\Phi,1}^{\mu,\mu}$ are not equivalent on $E.$
In particular the identity mapping from $E$ into $E$ can not be extended to an isomorphism between
$L_{\Phi,\alpha}(l_\infty,\mu,\nu)$ and $L_{\Phi,\alpha}(l_\infty,\mu,\mu).$

At the same time by following theorem the Orlicz spaces $L_{\Phi,\alpha}(l_\infty,\mu,\nu)$
and $L_{\Phi,\alpha}(l_\infty,\mu,\mu)$ are isometrically isomorphic.

\textbf{Theorem 4.} \emph{Let the    $N$-function $\Phi$ satisfy the
$(\delta_2,\Delta_2)$-condition, $\alpha\in [0,1].$ Then the Banach space $L_{\Phi,\alpha}(M,\varphi,\tau)$
is isometrically isomorphic to the Banach space}
$L_{\Phi}(M,\tau)=L_{\Phi,1}(M,\tau,\tau)$.

\emph{Proof.} For every $x\in
m_{\Phi,\alpha}^{\varphi,\tau}$ we have
$$
U(x)= (\Phi^{-1}(h))^{\alpha} x (\Phi^{-1}(h))^{1-\alpha}
\in L_\Phi(M,\tau).
$$
Therefore from definitions (1) and (4) of the norms we obtain
$$
\|x\|_{\Phi,\alpha}^{\varphi,\tau}=\|(\Phi^{-1}(h))^\alpha x
(\Phi^{-1}(h))^{1-\alpha}\|_\Phi.
$$
This means that the map $U$ defined as
$$
m_{\Phi,\alpha}^{\varphi,\tau}\ni
x\longrightarrow^{\!\!\!\!\!\!\!\!\!U} \  (\Phi^{-1}(h))^{\alpha} x (\Phi^{-1}(h))^{1-\alpha}\in
L_{\Phi}(M,\tau)\eqno(8)
$$
is a linear isometry. Let us show that the    $U(m_{\Phi, \alpha}^{\varphi, \tau})=(\Phi^{-1}(h))^{\alpha}
m_{\Phi,\alpha}^{\varphi,\tau} (\Phi^{-1}(h))^{1-\alpha}$ is dense in
$L_{\Phi}(M,\tau)$.

Let $h=\int\limits_{0}^{\infty} \lambda de_{\lambda}$ and
$q_{n}=\int\limits_{\frac{1}{n}}^{n}
de_{\lambda}\,\,,(n=1,2,...)$. Consider the set
$$
\mathcal{F}=\bigcup \limits_{m,n=1}^{\infty} q_{m}m_{\tau}q_{n}.
$$

Since $q_n\leq q_{n+1},$ it follows that $\mathcal{F}$ is a linear subspace in $m_\tau$ and by (2)
 $\mathcal{F}\subset L_{\Phi}(M,\tau).$

First. Let us prove that $\mathcal{F}$ is dense in $L_{\Phi}(M,\tau)$.
From the
$(\delta_2,\Delta_2)$-condition it follows that for $y\in L_{\Psi}(M,\tau)$
(where $\Psi$ is the complementary $N$-function for $\Phi$) the functional
$f(x)=\tau(xy),$ $x\in L_{\Phi}(M,\tau),$ defines the general form of continuous linear functional on
 $L_{\Phi}(M,\tau)$.

Let $y\in L_{\Psi}(M,\tau)$ and suppose that
$f(q_{m}xq_{n})=\tau((q_{m}xq_{n})y)=0$  for all $x\in
m_{\tau}$ and $m,n=1,2,... .$ In order to prove that $\mathcal{F}$ is dense in
$L_\Phi(M, \tau)$ it is sufficient to show that
$y=0.$

From the tracial property of $\tau$ we have that
$\tau(xq_{n}yq_{m})=0$ for all $x\in m_\tau.$
By proposition 2 $m_{\tau}$ is dense in
$L_{\Phi}(M,\tau)$ and hence  $q_{n}yq_{m}=0$ for all $m,n=1,2,...$.
Since $q_{n}\nearrow\textbf{1}$ as $n\rightarrow\infty,$
this implies that $y=0$.
Therefore $\mathcal{F}$ is dense
$L_{\Phi}(M,\tau)$.

Now let us show that $\mathcal{F}\subset U(m_{\Phi,\alpha}^{\varphi,\tau}).$ For this
it is sufficient to prove that given any  $x\in m_{\tau}$ and
$m,n=1,2,...,$ there exists $y\in m_{\Phi, \alpha}^{\varphi, \tau}$ such that
$q_{m}xq_{n}=U(y).$

Since the operators $(\Phi^{-1}(h))^{-\alpha}q_{m}$ and
$(\Phi^{-1}(h))^{\alpha-1}q_{n}$ belong to $M,$ the operator
$
y=U^{-1}(q_m xq_n)=(\Phi^{-1}(h))^{-\alpha}(q_{m}xq_{n})(\Phi^{-1}(h))^{\alpha-1}
$ also belongs to $M.$
From (2) and from  $\tau(|q_m x q_n|)<\infty$ we obtain that
$\tau(\Phi(|U(y)|))=\tau(\Phi(|q_m x q_n|))<\infty,$ i.e.
$y\in m_{\Phi, \alpha}^{\varphi, \tau}.$
 This implies that $\mathcal{F}\subset U(m_{\Phi,\alpha}^{\varphi,\tau}).$

Now since $m_{\Phi,\alpha}^{\varphi,\tau}$ is dense in
$(L_{\Phi,\alpha}(M,\varphi,\tau), \|\cdot\|_{\Phi,\alpha}^{\varphi,\tau})$ and
$U\left(m_{\Phi,\alpha}^{\varphi,\tau}\right)$ is dense in
$(L_\Phi(M,\tau), \|\cdot\|_{\Phi})$ the isometry
$U: m_{\Phi,\alpha}^{\varphi,\tau}\rightarrow L_\Phi(M,\tau)$  defined in (8)
can be uniquely extended to an isometric isomorphism between  $L_{\Phi,\alpha}(M,\varphi,\tau)$ and $L_{\Phi}(M,\tau).$
The proof is complete.\ $\Box$\\[-2mm]

Since every  faithful normal semi-finite trace $\tau_1$ on $M$ is a locally finite weight \cite{Tru3}
the theorem 4 implies the following

\textbf{Corollary 2.}\emph{ If $\tau_1$ and $\tau$  are faithful normal semi-finite traces on a von Neumann algebra $M,$
$\Phi$ is an $N$-function satisfying the $(\delta_2,\Delta_2)$-condition,
then the Orlicz spaces $L_\Phi(M,\tau_1)$ and $L_\Phi(M,\tau)$ are isometrically isomorphic.}

Theorem 4 and Corollary 2 together imply the following theorem

\textbf{Theorem 5.}\emph{ Let $\tau_1$ and $\tau$ be faithful normal traces on a von Neumann algebra  $M,$
and let
$\varphi_1, \varphi_2$ be faithful normal locally finite weights on $M.$ Suppose that $\Phi$ is an
$N$-function satisfying the
$(\delta_2,\Delta_2)$-condition, $\alpha,\beta\in[0,1].$ Then the Orlicz spaces
  $L_{\Phi,\alpha}(M,\varphi_1,\tau_1)$ and $L_{\Phi,\beta}(M,\varphi_2,\tau_2)$ are isometrically isomorphic.}

Theorem 4 implies also the following

\textbf{Corollary 3.} \emph{Let $\Phi$ be an  $N$-function satisfying the
$(\delta_2,\Delta_2)$-condition and let $\Psi$ be the complementary  $N$-function for  $\Phi$, and $\alpha,\beta\in[0,1].$
Then the dual space $(L_{\Phi,\alpha}(M,\varphi,\tau))^*$ for the Orlicz space
$L_{\Phi,\alpha}(M,\varphi,\tau)$ is isometrically isomorphic to the space $L_{\Psi}(M,\tau).$
If moreover $\Psi$ also satisfies the $(\delta_2,\Delta_2)$-condition then
$(L_{\Phi,\alpha}(M,\varphi,\tau))^*$ is isometrically isomorphic to  $L_{\Psi,\beta}(M,\varphi,\tau)$
and the Banach space
$L_{\Phi,\alpha}(M,\varphi,\tau)$ is reflexive.}

Now let us give a representation of the space  $L_{\Phi,\alpha}(M,\varphi,\tau)$
by locally measurable operators in the case where $\varphi$
is a regular locally finite weight, and the $N$-function $\Phi$ satisfies
$(\delta_2,\Delta_2)$-condition.

Consider the following subset in the algebra  $LS(M)$
of locally measurable operators affiliated with the von Neumann algebra $M$:
$$
\mathcal{L}_{\Phi,\alpha}(M,\varphi,\tau)=\{x\in LS(M):
O_{\Phi,\alpha}^{\varphi,\tau}(x)<\infty\},
$$
and for each $x\in \mathcal{L}_{\Phi,\alpha}(M,\varphi,\tau)$ put
$$
\|x\|_{\Phi,\alpha}^{\varphi,\tau}=\inf\left\{\lambda\geq0:
O_{\Phi,\alpha}^{\varphi,\tau}\left(\frac{x}{\lambda}\right)\leq1\right\}.
$$
It is clear that
$$
m_{\Phi,\alpha}^{\varphi,\tau}=M\bigcap
\mathcal{L}_{\Phi,\alpha}(M,\varphi,\tau).
$$
Repeating the proof of the Theorems 2 and 3 and of Corollary 1 we obtain that
$\mathcal{L}_{\Phi,\alpha}(M,\varphi,\tau)$ is a linear subspace of
 $LS(M)$ and that $\|\cdot\|_{\Phi,\alpha}^{\varphi,\tau}
$ is a norm on $\mathcal{L}_{\Phi,\alpha}(M,\varphi,\tau).$

\textbf{Theorem 6.}\emph{ Let $\varphi$ be a  regular locally finite normal weight on $M$ and suppose that $\Phi$
is an $N$-function satisfying the $(\delta_2,\Delta_2)$-condition and
  $\alpha\in [0,1].$ Then
 $(\mathcal{L}_{\Phi,\alpha}(M,\varphi,\tau),$ $\|\cdot\|_{\Phi,\alpha}^{\varphi,\tau})$
is a Banach space and $m_{\Phi, \alpha}^{\varphi, \tau}$ is dense in
 $(\mathcal{L}_{\Phi,\alpha}(M,\varphi,\tau),$ $\|\cdot\|_{\Phi,\alpha}^{\varphi,\tau}).$}

\emph{Proof.} Let  $h$ be the Radon-Nikodym derivative of the weight
$\varphi$ with respect to the trace $\tau$, and
suppose that $h=\int\limits_0^\infty\lambda de_\lambda$ is the spectral resolution of the operator $h$.
From Theorem 1 it follows that the operators $h$ and $h^{-1}$
are locally measurable. Therefore the operators
$\Phi^{-1}(h)=\int\limits_0^\infty\Phi^{-1}(\lambda)de_\lambda$ and
$(\Phi^{-1}(h))^{-1}=\int\limits_0^\infty((\Phi^{-1}(\lambda)))^{-1}de_\lambda$ are also locally measurable.
There exists a linear isometry $U$ from
$m_{\Phi,\alpha}^{\varphi,\tau}$ into $L_\Phi(M,\tau)$ (see (8)), in particular $U$
is injective. Thus  there exists the converse map
 $U^{-1}$ for the map $U.$ Since $U(m_{\Phi,\alpha}^{\varphi,\tau})$
is dense in $L_\Phi(M,\tau)$ (see the proof of theorem 4), the converse map
 $U^{-1}$ can be extended to a linear isometry from $L_{\Phi}(M,\tau)$
to the closure of $m_{\Phi,\alpha}^{\varphi,\tau}$ in
$(\mathcal{L}_{\Phi,\alpha}(M,\varphi,\tau), \|\cdot\|_{\Phi,\alpha}^{\varphi,\tau}).$

Similar to the proof of Theorem 4 we obtain that
 $m_{\Phi,\alpha}^{\varphi,\tau}$ is dense in $(\mathcal{L}_{\Phi,\alpha}(M,\varphi,\tau),$
$\|\cdot\|_{\Phi,\alpha}^{\varphi,\tau}).$
Therefore $U^{-1}$ can be extended to a linear isometry from
$(L_{\Phi}(M,\tau),$ $ \|\cdot\|_{\Phi})$ onto
$(\mathcal{L}_{\Phi,\alpha}(M,\varphi,\tau),$ $ \|\cdot\|_{\Phi,\alpha}^{\varphi,\tau}),$
i.e.
$(\mathcal{L}_{\Phi,\alpha}(M,\varphi,\tau),$ $ \|\cdot\|_{\Phi,\alpha}^{\varphi,\tau}),$ is a Banach space.
The proof is complete.\ $\Box$\\[-2mm]

Theorem 6 implies that in the case where $h$ and $h^{-1}$
are locally measurable operators and the
$N$-function $\Phi$ satisfies $(\delta_2,\Delta_2)$-condition, the Orlicz space $L_{\Phi,
\alpha}(M,\varphi,\tau)$ can be described by locally measurable operators in the following form
$$
L_{\Phi,\alpha}(M,\varphi,\tau)=\mathcal{L}_{\Phi,\alpha}(M,\varphi,\tau)=(\Phi^{-1}(h))^{-\alpha}
L_\Phi(M,\tau)(\Phi^{-1}(h))^{\alpha-1}\subset LS(M).
$$

\begin{center}
\textbf{Acknowledgements}
\end{center}

\emph{The final version of this work was done within the framework of the Associateship Scheme of
the Abdus Salam International Centre for Theoretical Physics (ICTP), Triest. Italy. The
first author thanks ICTP for providing financial support and all facilities during his stay in ICTP
(July-August, 2011). This work is supported in part  by the DFG AL 214 136-1 project (Germany).
The first and the third authors would like to thank the Institute of Applied Mathematics of the
Bonn University for hospitality (April-May, 2011).}

\end{document}